\documentclass[a4paper, 14pt, noxcolor,draft]{article}
\usepackage{maple2e}
 
\DefineParaStyle{Maple Output}
\DefineParaStyle{Maple Output}
\DefineCharStyle{2D Math}
\DefineCharStyle{2D Output}

\usepackage[a4paper]{hyperref}

\usepackage{fancybox}
\usepackage{xy}
\usepackage{epsfig}
\usepackage[T1]{fontenc} 
\usepackage[latin1]{inputenc}
\usepackage{amsfonts,amssymb,amsmath}
\usepackage{graphics} 
\usepackage{multido, xcolor}

\newtheorem{theo}{\sc Theorem}[section]

\newtheorem{rema}[theo]{\sc Remark}

\def\qed{ \ \hfil$\square$}

\makeindex

\def \F {{\mathbb F}}

\newcommand{\Gl}{{\rm Gl}}
\newcommand{\GL}{{\rm GL}}
\newcommand{\GSp}{{\rm GSp}}

\newcommand{\SL}{{\rm SL}}

\newcommand{\Q}{{\mathbb Q}}

\newcommand{\Z}{{\mathbb Z}}

\font\teneusm=eusm10 \font\seveneusm=eusm7 
\font\fiveeusm=eusm5 
\newfam\eusmfam 
\textfont\eusmfam=\teneusm 
\scriptfont\eusmfam=\seveneusm 
\scriptscriptfont\eusmfam=\fiveeusm

\def\mat #1,#2,#3,#4,{\left({#1\atop #3}{#2\atop #4}\right)}
\def\bra#1,{{\left\lbrace {#1}\right\rbrace}}

\def \diag {{\rm diag}}

\def \Sp{{\rm Sp}}
\def\tm{\tilde m}

\def\l1{\langle}

\newcommand{\B}{\left(\begin{array}{cc}}
\newcommand{\E}{\end{array}\right)}

\def \fns{{${}^{*)}$}}

\newcommand{\comm}[1]
{\fns\marginpar{$\boxed
{\hskip-6pt
{\small {\sf 
\begin{tabular} {l}
 #1
\end{tabular}
}
}
}
$
}
}
\def \?  {\comm{check ?}}

\usepackage{euscript}
\let\scr=\EuScript

\let\mathcal=\scr           
\def\ang#1,{{\left\langle {#1}\right\rangle}} 
\newcommand{\ds}{\displaystyle}

\def\qed{ \ \hfil$\square$}

\def \F {{\mathbb F}}

\font\teneusm=eusm10 \font\seveneusm=eusm7 
\font\fiveeusm=eusm5 
\newfam\eusmfam 
\textfont\eusmfam=\teneusm 
\scriptfont\eusmfam=\seveneusm 
\scriptscriptfont\eusmfam=\fiveeusm

\font\tengothic=eufm10
\font\sevengothic=eufm7
\font\fivegothic=eufm5
\newfam\Gothic
\textfont\Gothic=\tengothic
\scriptfont\Gothic=\sevengothic
\scriptscriptfont\Gothic=\fivegothic

\def\mat #1,#2,#3,#4,{\left({#1\atop #3}{#2\atop #4}\right)}
\def\bra#1,{{\left\lbrace {#1}\right\rbrace}}

\def \diag {{\rm diag}}

\def \Sp{{\rm Sp}}

\def\l1{\langle}

\let\scr=\EuScript

\let\mathcal=\scr           

\def\vin{{ {\tiny \mid }  
\kern-7.29pt 
\bigcup }}

\def\ang#1,{{\left\langle {#1}\right\rangle}} 

\def \cds{{\cdot{\dots}\cdot}}

\newcommand{\cross}{\times}
\normalsize

\newcounter{ncours}{\setcounter{ncours} {1}}

\def\diag {\mathop{\rm diag}\nolimits}
\def\SL {\mathop{\rm SL}\nolimits}
\def\Sp {\mathop{\rm Sp}\nolimits}

\def\qed{\quad\hbox{\hskip 1pt\vrule width 4pt height 6pt
          depth 1.5pt\hskip 1pt}}

\input amssym.def
\input amssym.tex

\font\teneusm=eusm10 \font\seveneusm=eusm7 
\font\fiveeusm=eusm5 
\newfam\eusmfam 
\textfont\eusmfam=\teneusm 
\scriptfont\eusmfam=\seveneusm 
\scriptscriptfont\eusmfam=\fiveeusm 
\def\scr#1{{\fam\eusmfam\relax#1}}

\def\finishproclaim{\par\rm
    \ifdim\lastskip<\medskipamount\removelastskip
     \penalty55\medskip\fi}

\def\proofof#1:{\par\medskip
   \noindent{\it Proof of {\rm #1}}}

\def\Ref[#1]{\par\smallskip\hang\indent%
  \llap{\hbox to\parindent{[#1]\hfil\enspace}}%
     \ignorespaces}
\def\Item#1{\par\smallskip
  \hang\indent\llap{\hbox to\parindent
     {#1\hfill\enspace}}\ignorespaces}
\def\ItemItem#1{\par\indent\hangindent2\parindent
\hbox to\parindent{#1\hfill\enspace}\ignorespaces}

\def\arrowsim{\smash{\mathop{\longrightarrow}
 \limits^{\lower1.5pt \hbox{$\scriptstyle\sim$}}}}

 \def\diag{{\rm diag}}

 \def\SL{{\rm SL}}

\def\GL{{\rm GL}}

\begin{document}

\title{
Explicit Shimura's conjecture for $\Sp_3$ on a computer
}
\author{ Alexei  Panchishkin, Kirill Vankov\\
{\tt http://www-fourier.ujf-grenoble.fr/\~{}panchish}
\\
\small e-mail : panchish$@$mozart.ujf-grenoble.fr, 
FAX:  33 (0) 4 76 51 44 78}

\date{}

\maketitle

\noindent
\begin{abstract}
We   compute by a different  method 
the generating series in 
Shimura's conjecture for $Sp_3$, proved by Andrianov in 1967.
We develop formulas for the Satake spherical maps for $\Sp_n$  and  $\Gl_n$.


\end{abstract}
\tableofcontents
\thispagestyle{empty}

\section{A formula for the generating series of Hecke operators of $Sp_3$}
A classical method to produce $L$-functions for an algebraic group $G$ over $\Q$ 
uses the generating series 
$$
\sum_{n=1}^{\infty} \lambda_f(n)n^{-s}=\prod_{p \ {\rm primes}}\sum_{\delta=0}^{\infty} \lambda_f(p^\delta )p^{-\delta s},
$$
of the eigenvalues of Hecke operators
 on an automorphic form $f$ on $G$.
We study the generating series of Hecke operators $\mathbf{T}(n)$ for  the symplectic group $Sp_g$, when
$g=3$, and $\lambda_f(n)=\lambda_f(\mathbf{T}(n))$.

Let  $\Gamma=\Sp_g(\Z)\subset \SL_{2g}(\Z)$ be the Siegel modular group of genus $g$, and  
 $[\mathbf{p}]_g=p\mathbf{I}_{2g}=\mathbf{T}(\underbrace{p, \cdots, p}_{2g})$  be the  scalar Hecke operator for $Sp_g$. 
According to Hecke and Shimura, 
\begin{align*}&
D_{p}(X)=\sum_{\delta=0}^{\infty} \mathbf{T}(p^\delta )X^\delta
\\ &
=\begin{cases}
\ds \frac 1 {1-\mathbf{T}(p)X+p[\mathbf{p}]_1X^2},  
\hskip 1.4cm  \mbox{ if } g=1 \\ 
\hskip 5cm
 \mbox{ (see \cite{Hecke}, and \cite{Shi71}, Theorem 3.21),}
\\ & \\
\ds \frac{1-p^2[\mathbf{p}]_2X^2}
 {
1-\mathbf{T}(p)X+\{p\mathbf{T}_1(p^2)+p(p^2+1)[\mathbf{p}]_2\}X^2
-p^3[\mathbf{p}]_2\mathbf{T}(p)X^3+p^6[\mathbf{p}]_2^2X^4
}&  \\ \hskip5cm\mbox{ if } g=2 \mbox{ (see \cite{Sh}, Theorem 2),}  
\end{cases}
\end{align*}
where 
$\mathbf{T}(p)$, $\mathbf{T}_i(p^2)$ ($i=1, \cdots, g$) are the $g+1$ generators of the corresponding Hecke ring over $\Z$
 for  the symplectic group $Sp_g$, in particular, $\mathbf{T}_g(p^2)=[\mathbf{p}]_g$.

The case $g=3$ was treated for the first time by Andrianov in \cite{An67}.
In the present paper we find a different  method to compute
 the polynomials $E(X)$  
and $F(X)$ in $X$ of degree 6 and 8 
such that  $D_p(X)=E(X)/F(X)$,  through the generators  $\mathbf{T}(p), \mathbf{T}_1(p^2), \mathbf{T}_2(p^2), 
[\mathbf{p}]_3$ of the Hecke ring (Theorems \ref{ThP} and \ref{th:ThQ}).  
Actually this series was computed
for $g=3$ by hand in Andrianov's original paper \cite{An67} using the multiplication rules for Hecke operators.
Our computation works also for higher genus (see \cite{VaSp4}).

The explicit knowledge of the sum of the generating series of Hecke operators
$$
D_{p}(X)=\sum_{\delta=0}^{\infty} \mathbf{T}(p^\delta )X^\delta=E(X)/F(X)
$$
gives a relation  between the Hecke eigenvalues and the Fourier coefficients of a Hecke  eigenform $f$.
This link is needed for constructing an analytic continuation of $L$-function on $Sp_g$,
which was
done by A.N.Andrianov in the case $g=2$, see \cite{An74}.
An approach for constructing an analytic continuation of the spinor $L$-function on $Sp_3$
was indicated at the talk of the first author, see \cite{PaGRFA}.

Our result is based on  the use of the Satake spherical map $\Omega$ (see Section \ref{sec:omega} for details). We obtain the 
following  formula for the polynomial $P({X})=\Omega(E({X}))$: 

\begin{equation}
\begin{split}
P({X}) = P_3(x_0,\, x_1&,\, x_2,\, x_3,\, {X})=\\
 =1&-\left({\displaystyle\frac{sym_{2,1,1}}{p}}+{\displaystyle\frac{(p^2+p+1)\, sym_{1,1,1}}{p^2}}+{\displaystyle\frac{sym_{1,1,0}}{p  }}\right) \, x_0^2\, {X}^2\\
   &+{\displaystyle\frac{p+1}{p^2}}\left(sym_{2,2,2}+sym_{2,2,1}+sym_{2,1,1}+sym_{1,1,1}\right) \, x_0^3\, {X}^3\\ 
   &-\left({\displaystyle\frac{sym_{3,2,2}}{p^2}}+{\displaystyle\frac{(p^2+p+1)\, sym_{2,2,2}}{p^3}}+{\displaystyle\frac{sym_{2,2,1}}{p^2}}\right) \, x_0^4\, {X}^4\\
   &+{\displaystyle\frac{sym_{3,3,3}}{p^3}}\, x_0^6\,{X}^6 \, .
\end{split}
\end{equation}
(see also \cite{PaVa}).
The case $g=4$ was treated by K.Vankov in  \cite{VaSp4}.

\smallskip
In this formula the notation $sym_{i_1,i_2,i_3}$ represents the symmetric polynomial of three variables 
$x_1$, $x_2$ and $x_3$ constructed in the following way:
$$
sym_{i_1,i_2,i_3} = \sum_{\sigma\in S_n / \,{\rm Stab}(x_1^{i_1}x_2^{i_2}x_3^{i_3})} \sigma (x_1^{i_1}x_2^{i_2}x_3^{i_3}),
$$
where the summation of permuted monomials $x_1^{i_1}x_2^{i_2}x_3^{i_3}$ is normalized using the stabilizer
${\rm Stab}(x_1^{i_1}x_2^{i_2}x_3^{i_3})$ so that all coefficients are equal to 1 and $i_1\ge i_2\ge i_3\ge 0$.  
The total degree of the polynomial is $i_1+i_2+i_3$.  Here $S_n = S_3$ is the symmetric group that acts naturally on polynomials in $n$ variables, where $n=3$ in our case. For example
\begin{equation}\nonumber
\begin{split}
sym_{0,0,0}&=1\\
sym_{1,0,0}&=x_1+x_2+x_3\\
sym_{1,1,0}&=x_1x_2+x_1x_3+x_2x_3\\
sym_{1,1,1}&=x_1x_2x_3\\
sym_{4,3,2}&=x_1^4x_2^3x_3^2+x_1^4x_2^2x_3^3+x_1^3x_2^4x_3^2+x_1^3x_2^2x_3^4+x_1^2x_2^4x_3^3+x_1^2x_2^3x_3^4 \, .
\end{split}
\end{equation}
Many computations presented in this article were performed using Maple 9.50 (IBM INTEL NT).  
Symmetric polynomials $sym_{i_1i_2i_3}$ (up to total degree 9) were computed using 
the coefficient of $t$ of the generating function 
\begin{equation}\nonumber
\begin{split}
&\prod_{\sigma\in S_3}(1+tx_{\sigma(1)}^{i_1}x_{\sigma(2)}^{i_2}x_{\sigma(3)}^{i_3})=(1+tx_1^{i_1}x_2^{i_2}x_3^{i_3})(1+tx_1^{i_1}x_2^{i_3}x_3^{i_2}) \\
&\qquad \cross (1+tx_1^{i_2}x_2^{i_1}x_3^{i_3})(1+tx_1^{i_2}x_2^{i_3}x_3^{i_1})(1+tx_1^{i_3}x_2^{i_1}x_3^{i_2})(1+tx_1^{i_3}x_2^{i_2}x_3^{i_1}) \, ,
\end{split}
\end{equation}
where $i_1=0,\dots ,6$, $i_2=0,\dots ,i_1$ and $i_3=0,\dots ,i_2$. Then it is normalized by dividing out 
its leading coefficient.

Let us state our result directly in terms of the Hecke operators for the symplectic group $\Sp_n$, 
defined at p. 142 of \cite{An87}.
Consider the group of positive symplectic similitudes
\begin{align}
& \mathrm{S}=\mathrm{S}^n=\GSp_n^+(\Q)=\{M\in \mathrm{M}_{2n}(\Q)\ |\ {}^tMJ_nM=\mu(M)J_n, \mu(M)>0\} \ ,
\\ & 
\nonumber 
\mbox{where } J_n = \left( \begin{array}{cc}\mathbf{0}_n& \mathbf{I}_n \\ -\mathbf{I}_n& \mathbf{0}_n\end{array}\right) . 
\end{align}
For the  Siegel  modular group 
$\Gamma=\Sp_n(\Z)$ consider the double cosets
\begin{align}&
(M)=\Gamma M\Gamma \subset 
\mathrm{S}, 
\end{align}
and the Hecke operators
\begin{align}&
\mathbf{T}(a)= \sum_{M\in \mathrm{SD}_n(a)}(M), 
\end{align}
where $M$ runs through the following integral matrices 
\begin{align}&
\mathrm{SD}_n(a)=\{\diag(d_1, \cdots, d_n;e_1, \cdots, e_n)\ |\ d_i| d_{i+1}, d_n| e_{n}, e_{i+1}|e_i, d_ie_i=a\}. 
\end{align}
Let us use the notation
\begin{align}&
\mathbf{T}(d_1, \cdots, d_n;e_1, \cdots, e_n)=(\diag(d_1, \cdots, d_n;e_1, \cdots, e_n)).
\end{align}
In particular we have the operators
(see p. 149 of \cite{An87}):
\begin{align}&
\mathbf{T}(p)= \mathbf{T}(\underbrace{1, \cdots, 1}_{n},\underbrace{p, \cdots, p}_{n}),\\ &
\mathbf{T}_i(p^2)= \mathbf{T}(\underbrace{1, \cdots, 1}_{n-i},\underbrace{p, \cdots, p}_{i},\underbrace{p^2, \cdots, p^2}_{n-i},\underbrace{p, \cdots, p}_{i}),
\mbox{ for } i=1, 2, \cdots, n.
\end{align}
Then their images by the spherical map $\Omega$ are given at p.159 of  \cite{An87}:
\begin{align}&
\Omega(\mathbf{T}(p))= x_0\prod_{i=0}^n(1+x_i)=
\sum_{j=0}^nx_0s_j(x_1, x_2, \cdots, x_n),\\ &
\Omega(\mathbf{T}_i(p^2))= \sum_{a+b\le n, a\ge i}p^{b(a+b+1)}\mathrm{sm}_p(a-i, a)x_0^2\omega(\pi_{a, b}).
\end{align}
Here
$$s_i(x_1,\dots ,x_n) = \sum_{1\le \alpha_1 < \cdots < \alpha_i \le n} x_{\alpha_1} \cdots x_{\alpha_i}$$
is the $i$th elementary symmetric polynomial (different then previously defined \newline
$sym_{i_1,i_2,i_3}$),
$\pi_{a, b}=\left(\left(\begin{array}{ccc}\mathbf{I}_{n-a-b} & &\\ & p\mathbf{I}_{a}&\\ & &p^2\mathbf{I}_{b}\end{array}\right)\right)$ is a Hecke operator for $\GL_n$, and
the coefficient $\mathrm{sm}_p(r, a)$ denotes the number of symmetric matrices of rank $r$ and order $a$ over the field $\F_p$.
This coefficient is evaluated at p.205 of  \cite{An87}:
\begin{align}\label{phi}&
\mathrm{sm}_p(r, a)=\mathrm{sm}_p(r, r)\frac{\phi_a(p)}{\phi_{r}(p)\phi_{a-r}(p)}, 
\\ & \nonumber
\mbox{ with  } \phi_{r}(x)=(x-1)(x^2-1)\cds(x^r-1) \mbox{ for   } r\ge 1, 
\mbox{ and  } \phi_{0}(x)=1.
\end{align}
In particular, we have in the case $n=3$ that 
\begin{align}\label{Tp}
\nonumber
&\Omega(\mathbf{T}(p)) 
  = \mathit{x}_0\,(1 + {\mathit{sym}_{1, \,0, \,0}} + {\mathit{sym}_{1, \,1, \,0}} + {\mathit{sym}_{1, \,1, \,1}}) \ ,
\\   
& \Omega(\mathbf{T}_1(p^2))  = {\displaystyle \frac {\mathit{x}_0^{2}\,(p^{2} - 1)}{p^{3}}\Big( {\mathit{sym}_{2, \,1, \,1}} 
+ {\mathit{sym}_{1, \,1, \,0}} \Big)}
\\ \nonumber &
  + {\displaystyle \frac {\mathit{x}_0^{2}}{p}}\Big( {\mathit{sym}_{2, \,2, \,1}}
+ {\mathit{sym}_{2, \,1, \,0}} + {\mathit{sym}_{1, \,0, \,0}} \Big)
\\ \nonumber 
 &\quad 
+ {\displaystyle \frac {\mathit{x}_0^{2}
(p-1)(3p^2+2 p+1)
}
{p^{4}}}
{\mathit{sym}_{1, \,1, \,1}} \ ,
\\ \nonumber
&\Omega(\mathbf{T}_2(p^2))
  = p^{0}\mathrm{sm}_p(0, 2)x_0^2\omega(\pi_{2, 0}) + p^{4}\mathrm{sm}_p(0, 2)x_0^2\omega(\pi_{2, 1}) 
\\ \nonumber
&
+ p^{0}\mathrm{sm}_p(1, 3)x_0^2\omega(\pi_{3, 0})
\\ \nonumber
 &\quad = x_0^2\omega (t(1,p, p))+p^{4}x_0^2\omega (t(p, p, p^2))+\mathrm{sm}_p(1, 3)x_0^2\omega (t(p, p, p))
\\ \nonumber
 & \quad = {\displaystyle \frac{\mathit{x}_0^{2}}{p^3}} \left( \mathit{sym}_{1, \,1, \,0} + \mathit{sym}_{2, \,1, \,1} \right)
  + {\displaystyle \frac{\mathit{x}_0^{2}\,(p-1)(p^2+p+1)}{p^{6}}} \mathit{sym}_{1, \,1, \,1} \ ,  
\\ \nonumber
\\ \nonumber
&\Omega(\mathbf{T}_3(p^2))
  = \Omega([\mathbf{p}]_3)
   = p^{0}\mathrm{sm}_p(0, 3)x_0^2\omega(\pi_{3, 0}) 
   = \frac{x_0^2x_1x_2x_3}{p^{6}}
   = {\displaystyle \frac{\mathit{x}_0^{2}}{p^{6}}} sym_{1,1,1} \ , 
\end{align}
because of the equality (\ref{phi}) with $a=3, r=1$ implying $\mathrm{sm}_p(1, 3)= (p-1)(p^2+p+1)$.

Consider the polynomial $Q_3({X})$ defining the spinor zeta function
$Z(s)$ of genus three:
\begin{align}
\nonumber
Q_3({X}) &= Q_3(x_0, x_1, x_2, x_3, {X}) 
\\ 
&= (1-x_0{X})(1-x_0x_1{X})(1-x_0x_2{X})(1-x_0x_3{X})
\\ \nonumber
& \cross (1-x_0x_1x_2{X})(1-x_0x_1x_3{X})
(1-x_0x_2x_3{X})(1-x_0x_1x_2x_3{X}) .
\end{align}
Following the proof at  p.159 of  \cite{An87}, there exist Hecke operators  
$$
\mathbf{q}_j\in 
\Q[ \mathbf{T}(p), \mathbf{T}_1(p^2), \cdots, \mathbf{T}_{n}(p^2)]\mbox{ such that}
$$ 
\begin{align}&
\sum_{j=0}^{2^n}\Omega(\mathbf{q}_j){X}^j = Q_n({X})
\\ \nonumber &
=
(1-x_0{X})(1-x_0x_1{X})(1-x_0x_2{X})\cds (1-x_0x_1x_2\cds x_n{X}) .
\end{align}
Let us consider the series 
$D({X})  = \sum_{\delta=0}^{\infty} \mathbf{T}(p^\delta ){X}^\delta  \in {\mathcal L}_\Z[\![{X}]\!]$ and 
the polynomial $F({X})=\sum_{j=0}^{2^n}\mathbf{q}_j{X}^j$ over the Hecke ring ${\mathcal L}_\Z$.

It was established by A.N.Andrianov, 
that there exist polynomials  
$$
E({X}) 
\in \Q[ \mathbf{T}(p), \mathbf{T}_1(p^2), \cdots, \mathbf{T}_{n}(p^2), {X}] \mbox{ such that }
$$
\begin{align}&
D({X})  = \sum_{\delta=0}^{\infty} \mathbf{T}(p^\delta ){X}^\delta  =
\frac{E({X})}{F({X})},
\end{align}
with the above polynomial $F({X})=\sum_{j=0}^{2^n}\mathbf{q}_j{X}^j$ of degree $2^n$, and  such that
 $\ds E({X})=\sum_{j=0}^{2^n-2}\mathbf{u}_j{X}^j$ is a polynomial
 of degree $2^n-2$ 
with the leading term 
$$
(-1)^{n-1}p^{{n(n+1)}2^{n-2}-n^2}[\mathbf{p}]^{2^{n-1}-1} {X}^{2^n-2}
$$
(as stated in  Theorem 6 at p. 451 of \cite{An70} and at p.61 of \S 1.3, \cite{An74}).
In the following theorem we denote by $[\mathbf{p}]_n=(p\mathbf{I}_{2n})=\mathbf{T}_n(p^2)$  the  element (3.4.48) in \cite{An70}), so that 
$\Omega([\mathbf{p}]_n)= p^{-n(n+1)/2}x_0^2x_1\cds x_n$.
\begin{theo}[see also  \cite{An67}]\label{ThP} If $n=3$, there is the following explicit polynomial presentation: 
\begin{align}&
D({X})  = \sum_{\delta=0}^{\infty} \mathbf{T}(p^\delta ){X}^\delta 
=\frac{E({X})}{F({X})},
\end{align}
where
\begin{align}\label{P3}& \nonumber
E({X}) 
\\  \quad &
 =1-
p^2\left( \mathbf{T}_2(p^2)+
(p^2-p+1)(p^2+p+1)
[\mathbf{p}]_3\right){X}^2
+(p+1)p^4\mathbf{T}(p)[\mathbf{p}]_3{X}^3
\\ \nonumber \quad
      &
-p^7[\mathbf{p}]_3\left( \mathbf{T}_2(p^2)+
(p^2-p+1)(p^2+p+1)
[\mathbf{p}]_3\right) {X}^4 
   +p^{15}[\mathbf{p}]_3^3\,{X}^6 \, \in {\mathcal L}_\Z[{X}].
\end{align}


\end{theo}
\begin{rema} It was pointed out to the first author by S.Boecherer, that
difficulties in the problem of analytic continuation of the spinor $L$-function of genus 3
could come from
the polynomial $E({X})$.
Indeed, this is clearly indicated  by Kurokawa's paper \cite{Ku88}, Theorem 2
in the case of the Siegel-Eisenstein series of genus 3.
A similar polynomial is discussed in the paper of Maass \cite{Maa76}
using densities.
\end{rema}
\begin{rema} It seems that the polynomial $E({X})$ 
does not depend on $\mathbf{T}_1(p^2)$ in general.
Also,  we conjecture that the coefficients of $E(X)$
at $X$ and at $X^{2^n-3}$ are always equal to zero.
\end{rema}

\noindent{\it Proof} 
of Theorem  \ref{ThP} is completed in Section \ref{sec:omega}.

\section{Explicit form of Shimura's conjecture for $Sp_3$}
We derive a different  method to compute
the generating series in Shimura's conjecture for $Sp_3$
first computed  by A.N.Andrianov 
in \cite{An67} (see  also \cite{An68} and \cite{An69}).
Shimura's conjecture was stated in  \cite{Sh}, at p.825 as follows:
\begin{quote}
``In general, it is plausible that $D_p(X)=E(X)/F(X)$ with polynomials $E(X)$  
and $F(X)$ in $X$ with integral coefficients of degree $2^n-2$ and $2^n$, respectively'' 
\end{quote}
({\it i.e. with coefficients in } ${\mathcal L}_\Z=\Z[ \mathbf{T}(p), \mathbf{T}_1(p^2), \cdots, \mathbf{T}_{n}(p^2)]$).
\begin{theo}[see also  \cite{An67}]\label{th:ThQ} If $n=3$, one has the following explicit polynomial presentation: 
\begin{align}\label{ThQ}&
D(X)  = \sum_{\delta=0}^{\infty} \mathbf{T}(p^\delta ){X}^\delta 
=\frac{E(X)}{F(X)},
\end{align}
where $E(X)$ is given by (\ref{P3}), and 
\begin{align}&
F(X)
= 1 - \mathbf{T}(p)\,{X} 
\\ & \nonumber
+p\big(\,
\mathbf{T}_1(p^2) 
+ (p^{2} + 1)\,\mathbf{T}_2(p^2) 
+ 
(p^2+1)^2
[\mathbf{p}]_3\big)\,{X}^{2} \\ & 
\nonumber \mbox{}  
- p^{3}\mathbf{T}(p)\,\big( \mathbf{T}_2(p^2)
+ [\mathbf{p}]_3\big)\,{X}^{3} 
\\ & \nonumber
+ p^6\big(
 - 2\,p\,\mathbf{T}_1(p^2)\,[\mathbf{p}]_3 + \,\mathbf{T}_2(p^2)
\mbox{}  - 2( p -1)\,\mathbf{T}_2(p^2)\,[\mathbf{p}]_3
 \\ &\nonumber
-(p^2+2p-1)(p^2-p+1)(p^2+p+1))
[\mathbf{p}]_3^{2} +
[\mathbf{p}]_3\,
\mathbf{T}(p)^{2}\big){X}^{4} \\ &\nonumber
\mbox{} 
- \,p^{9}\,[\mathbf{p}]_3^{}\,\mathbf{T}(p)\,\big ( 
\mathbf{T}_2(p^2) +[\mathbf{p}]_3\big )\,
{X}^{5}   \\ &\nonumber
+p^{13}\,[\mathbf{p}]_3^{2}\,\big(\,\mathbf{T}_1(p^2) + (
p^{2} + 1)\,\mathbf{T}_2(p^2) + 
(1+p^2)^2
[\mathbf{p}]_3\big)\,{X}^{6} \\ &\nonumber
\mbox{} - \,p^{18}\,[\mathbf{p}]_3^{3}\,\mathbf{T}(p)
\,{X}^{7} +\,p^{24}\,[\mathbf{p}]_3^{4}\,{X}^{8}
\in {\mathcal L}_\Z[{X}].
\end{align}

\end{theo}

\begin{rema}
(a)
In the case $n=2$ Shimura proved (\cite{Sh},  Theorem 2)  that
\begin{align*}&
\sum_{\delta=0}^{\infty} \mathbf{T}(p^\delta )X^\delta = (1-p^2
[\mathbf{p}]_2X^2)\times
\\ & \qquad
[1-\mathbf{T}(p)X+\{p\mathbf{T}_1(p^2)+p(p^2+1)[\mathbf{p}]_2\}X^2
-p^3[\mathbf{p}]_2\mathbf{T}(p)X^3+p^6
[\mathbf{p}]_2^2X^4
]^{-1}
\end{align*}

\noindent (b)
  For the group $G=\GL_n$ it was proved by T.Tamagawa  \cite{Tam}  that  for all $n$
$$
\sum_{\delta=0}^{\infty} t(p^\delta )X^\delta=[\sum_{i=0}^n(-1)^i p^{i(i-1)/2}\pi_i(p)X^i]^{-1}
$$
(see in \cite{Shi71},  Theorem 3.21).
\end{rema}

\begin{proof}
of Theorem  \ref{th:ThQ} follows the same lines as that of Theorem \ref{ThP}.
We compute an expression for $\Omega(F({X}))$ in terms of
$sym_{i_1,i_2,i_3}$:
\begin{align}& \nonumber
\Omega(F({X}))
=1 - \mathit{x}_0^{}\,({\mathit{sym}_{1, \,1, \,1}} + {\mathit{sym}_{1, \,1, \,0}} + {
\mathit{sym}_{1, \,0, \,0}} +  1)\,
{X} \\ &  \nonumber
\mbox{} + \mathit{x}_0^{2}\,(4\,{\mathit{sym}_{1, \,1, \,1}} + {
\mathit{sym}_{1, \,0, \,0}} + 2\,{\mathit{sym}_{2, \,1, \,1}} + 2
\,{\mathit{sym}_{1, \,1, \,0}} 
\\ &  \nonumber\hskip3cm
+ {\mathit{sym}_{2, \,1, \,0}} + {
\mathit{sym}_{2, \,2, \,1}})\,{X}^{2}  \\ &  \nonumber
-\mathit{x}_0^{3}({\mathit{sym}_{3, \,1, \,1}} + {\mathit{sym}_{1, 
\,1, \,0}} + 4\,{\mathit{sym}_{2, \,2, \,1}} + 4\,{\mathit{sym}_{
1, \,1, \,1}} + {\mathit{sym}_{2, \,1, \,0}} 
\\ &  \nonumber+ {\mathit{sym}_{2, 
\,2, \,0}} + 4\,{\mathit{sym}_{2, \,2, \,2}} 
\mbox{} + {\mathit{sym}_{3, \,2, \,2}} + {\mathit{sym}_{3, \,2, 
\,1}} + 4\,{\mathit{sym}_{2, \,1, \,1}}){X}^{3}\mbox{} 
\\ &  \nonumber
+ 
\mathit{x}_0^{4}({\mathit{sym}_{3, \,1, \,1}} + {\mathit{sym}_{1, \,1, \,1}}
 + {\mathit{sym}_{3, \,3, \,1}} \\ &  
\mbox{} + {\mathit{sym}_{4, \,2, \,2}} + 2\,{\mathit{sym}_{2, \,1
, \,1}} + 4\,{\mathit{sym}_{3, \,2, \,2}} + 2\,{\mathit{sym}_{3, 
\,2, \,1}} + {\mathit{sym}_{2, \,2, \,0}} + 8\,{\mathit{sym}_{2, 
\,2, \,2}} \\ &  \nonumber
\mbox{} + 2\,{\mathit{sym}_{3, \,3, \,2}} + {\mathit{sym}_{3, \,3
, \,3}} + 4\,{\mathit{sym}_{2, \,2, \,1}}){X}^{4}\mbox{} 
\\ &  \nonumber
- 
\mathit{x}_0^{5}({\mathit{sym}_{4, \,3, \,3}} + {\mathit{sym}_{4, \,3, \,2
}} + {\mathit{sym}_{2, \,2, \,1}} \\ &  \nonumber
\mbox{} + 4\,{\mathit{sym}_{2, \,2, \,2}} + 4\,{\mathit{sym}_{3, 
\,3, \,2}} + {\mathit{sym}_{3, \,3, \,1}} + 4\,{\mathit{sym}_{3, 
\,2, \,2}} 
\\ &  \nonumber\hskip3cm
+ 4\,{\mathit{sym}_{3, \,3, \,3}} + {\mathit{sym}_{4, 
\,2, \,2}} + {\mathit{sym}_{3, \,2, \,1}} 
){X}^{5}
 \\ &  \nonumber
 + \mathit{x}_0^{6}\,(2\,{\mathit{sym}_{3, \,3, \,2}}
 + {\mathit{sym}_{3, \,2, \,2}} + 2\,{\mathit{sym}_{4, \,3, \,3}}
\\ &  \nonumber\hskip3cm 
+ 4\,{\mathit{sym}_{3, \,3, \,3}} + {\mathit{sym}_{4, \,3, \,2}}
 + {\mathit{sym}_{4, \,4, \,3}})\,{X}^{6} \\ &  \nonumber
\mbox{} - \mathit{x}_0^{7}\,({\mathit{sym}_{4, \,3, \,3}} + {
\mathit{sym}_{3, \,3, \,3}} + {\mathit{sym}_{4, \,4, \,4}} + {
\mathit{sym}_{4, \,4, \,3}})\,{X}^{7} + \mathit{x}_0^{8}\,{\mathit{
sym}_{4, \,4, \,4}}\,{X}^{8}.
 \end{align}
Then we use the polynomial expressions (\ref{Tp}) 
for the generators of the Hecke ring.
Using these generators,  we may look for a solution 
in the following form:
\begin{align}&
\Omega(F({X}))
= 1 - \Omega(\mathbf{T}(p)){X}
 + 
\big(\mathit{K}_{T1p2}\Omega(\mathbf{T}_1(p^2)) + \mathit{K}_{T2p2}\Omega(\mathbf{T}_2(p^2))
 \\ & \nonumber
 + \mathit{K}_{T3p2}\Omega([\mathbf{p}]_3) + \mathit{K}_{TpTp}\Omega(\mathbf{T}(p))^{2}\big)\,{X}^{2} 
  \\ & \nonumber
+\big(\mathit{K}_{TpT1p2}\Omega(\mathbf{T}(p)\mathbf{T}_1(p^2)) 
+ 
\mathit{K}_{TpT2p2}\Omega(\mathbf{T}(p)\mathbf{T}_2(p^2)) 
 \\ & \nonumber
+ \mathit{K}_{TpT3p2}\,
\Omega(\mathbf{T}(p)[\mathbf{p}]_3) + \mathit{K}_{TpTpTp}\Omega(\mathbf{T}(p))^{3}\big)
{X}^{3}\mbox{}
 \\ & \nonumber
 + 
\big(\mathit{K}_{T1p2T1p2}\Omega(\mathbf{T}_1(p^2))^{2} 
+ \mathit{K}_{T1p2T2p2}\Omega(\mathbf{T}_1(p^2)\mathbf{T}_2(p^2)) 
\\ & \nonumber
+ 
\mathit{K}_{T1p2T3p2}\,\Omega(\mathbf{T}_1(p^2)[\mathbf{p}]_3)
\mbox{} + \mathit{K}_{T2p2T2p2}\Omega(\mathbf{T}_2(p^2))^{2}
\\ & \nonumber
 + 
\mathit{K}_{T2p2T3p2}\,\Omega(\mathbf{T}_2(p^2)[\mathbf{p}]_3) + \mathit{K}_{T3p2T3p2}\,
\Omega([\mathbf{p}]_3^{2}) \\ &\nonumber
\mbox{} + \mathit{K}_{T1p2TpTp}\,\Omega(\mathbf{T}_1(p^2)\mathbf{T}(p)^{2}) + 
\mathit{K}_{T2p2TpTp}\Omega(\mathbf{T}_2(p^2)\mathbf{T}(p)^{2}) 
\\ &\nonumber
+ \mathit{K}_{T3p2TpTp}\,\Omega([\mathbf{p}]_3\mathbf{T}(p)^{2}) 
\mbox{} + \mathit{K}_{TpTpTpTp}\,\Omega(\mathbf{T}(p))^{4}\big){X}^{4}\\
&\nonumber
+
\Omega(\mathbf{q}_5){X}^5+\Omega(\mathbf{q}_6){X}^6
+\Omega(\mathbf{q}_7){X}^7+\,p^{24}\,\Omega([\mathbf{p}]_3)^{4}{X}^8.  
\end{align} 

It is not too difficult to resolve the resulting equations in the indeterminate coefficients:\\
$
K_{TpT1p2} = 0, K_{TpTpTp} = 0, K_{TpT2p2} = -p^3, K_{TpT3p2} =
-p^3,
K_{T2p2} = p^3+p$, \\ $K_{T3p2} = 
p(1+p^2)^2,
K_{T1p2} = p, K_{TpTp} = 0$,
$
K_{T2p2TpTp}=0$, \\ 
$K_{T1p2T3p2}=-2p^7,
K_{T2p2T3p2}= -2p^7+2p^6, 
K_{T1p2T1p2}= 0,
K_{T1p2T2p2}= 0$, \\ 
$K_{T2p2T2p2}= p^6$, 
$K_{T1p2TpTp}= 0, 
K_{T3p2TpTp}=p^6, 
K_{TpTpTpTp}= 0$,\\ $
K_{T3p2T3p2}=
-p^6(p^2+2p-1)(p^2-p+1)(p^2+p+1)$.

Then we find the remaining coefficients $\mathbf{q}_5$, $\mathbf{q}_6$, $\mathbf{q}_7$ using the functional equation \cite{An87}, p.164
(3.3.79):   
$\mathbf{q}_{8-i}=(p^6[\mathbf{p}]_3)^{4-i}\mathbf{q}_i$ ($i=0, \cdots, 8$), 
compare with formulas in 
\cite{An67} and in 
\cite{Evd}. 
\end{proof}

\section{An identity involving $\omega (t(1,p^{\lambda_2}, p^{\lambda_3}))$}\label{sec:omega}
The theory of Hecke rings for the symplectic group is developed in \cite{Sh}, 
\cite{An87} and \cite{AnZh95} (Ch. 3).  
At page 150 of \cite{AnZh95} we have the following identity for the spherical map 
\begin{equation}
\begin{split}
R_n({X}) & = \sum_{\delta=0}^{\infty} \Omega (\mathbf{T}(p^\delta )){X}^\delta \\
 & = \sum_{\delta=0}^{\infty}\quad\sum_{0\le\delta_1\le\cdots\le\delta_n\le\delta}p^{n\delta_1+(n-1)\delta_2+\cdots+\delta_n}\omega(t(p^{\delta_1},\cdots,p^{\delta_n}))(x_0{X})^\delta \, .
\end{split}
\end{equation}

This identity for formal generating series of elements of Hecke ring allows 
to reduce computations in the local Hecke rins of the symplectic group 
to computations in polynomial rings by applying the spherical map $\Omega$ 
to elements $\mathbf{T}(p^\delta ) = \mathbf{T}^n(p^\delta ) \in 
L_{\Q}(\Gamma^n,S^n)
$ of Hecke ring for the symplectic 
group or spherical map $\omega$ to elements  
$t(p^{\delta_1},\cdots,p^{\delta_n}) = \sum_j a_j(\Lambda g_j) \in L_{\Q}(\Lambda ^n,G^n)$ 
of Hecke ring for the general linear group. 
Detailed definition of spherical maps 
for Hecke elements as well as definitions and a structure of left cosets $\Lambda g_j$ 
and Hecke pairs $(\Gamma^n,S^n)$ and $(\Lambda ^n,G^n)$ can be found in \cite{AnZh95},  chapter 3 paragraphs 
2 and 3, see also \cite{An87} chapter 3.  
For generating elements $\pi_i(p) = \pi_i^n(p) = (\diag(1,\dots ,1,p,\dots ,p))$ with 1 
on the diagonal listed $(n-i)$ times and then $p$ listed $i$ times ($1\le i \le n$) 
the images elements under the map $\omega = \omega_p^n$ are given by the formulas
\begin{equation}\nonumber
\omega (\pi_i^n(p))=p^{- \langle i \rangle }s_i(x_1,\dots ,x_n) \quad (1\le i \le n),
\end{equation}
where
\begin{equation}\nonumber
s_i(x_1,\dots ,x_n) = \sum_{1\le \alpha_1 < \cdots < \alpha_i \le n} x_{\alpha_1} \cdots x_{\alpha_i}
\end{equation}
is the $i$th elementary symmetric polynomial. 
For an arbitrary element $t$ the map is defined by
\begin{equation}\nonumber
\omega(t)=\sum_j a_j\omega ((\Lambda g_j)) \, .
\end{equation}

Examples of computation for cases $n=1$ and $n=2$ are given in \cite{Sh}, at page 824, 
and in the book \cite{AnZh95}, at page 150:
\begin{align}&
R_1({X})=[(1-x_0{X})(1-x_0x_1{X})]^{-1}\, ,\\
&
R_2({X})=\frac{ 1-p^{-1}x_0^2x_1x_2{X}^2}{(1-x_0{X})(1-x_0x_1{X})(1-x_0x_2{X})(1-x_0x_1x_2{X})} \, .
\end{align}

\begin{proof}
of Theorem  \ref{ThP}.
Let us consider $n=3$.  Acting analogously to the case $n=2$ let 
\begin{equation}
\begin{split}
\delta_2 &=\delta_1+\delta_1^\prime\\
\delta_3 &=\delta_1+\delta_2^\prime\\
\delta &=\delta_1+\delta^\prime\\
\delta^\prime &=\delta_2^\prime+\beta
\end{split}
\end{equation}
where $0\le\delta_1^\prime\le\delta_2^\prime\le\delta^\prime$, $\beta\ge0$.  Then
\begin{equation}\nonumber
\begin{split}
R_3({X})&=\sum_{\delta =0}^\infty \, \sum_{0\le \delta_1 \le \delta_2 \le \delta_3 \le \delta} p^{3\delta_1 +2\delta_2 +\delta_3} \omega (t(p^{\delta_1},p^{\delta_2},p^{\delta_3})) (x_0 {X})^\delta\\
&=\sum_{\delta_1\ge 0}\sum_{\beta\ge 0\atop 0\le\delta_1^\prime\le\delta_2^\prime} (x_0{X})^{\delta_1+\delta_2^\prime+\beta} \left(\frac{x_1x_2x_3}{p^6}\right)^{\delta_1} \omega(t(1,p^{\delta_1^\prime},p^{\delta_2^\prime})) p^{3\delta_1+2(\delta_1+\delta_1^\prime)+(\delta_1+\delta_2^\prime)}\\
&=\sum_{\delta_1\ge 0}\sum_{\beta\ge 0, 0\le\delta_1^\prime\le\delta_2^\prime} (x_0{X})^{\delta_1+\delta_2^\prime+\beta} \left(\frac{x_1x_2x_3}{p^6}\right)^{\delta_1} p^{6\delta_1+2\delta_1^\prime+\delta_2^\prime} \omega(t(1,p^{\delta_1^\prime},p^{\delta_2^\prime}))\\
&=\sum_{\delta_1\ge 0}\sum_{\beta\ge 0} (x_0{X}x_1x_2x_3)^{\delta_1} (x_0{X})^\beta \sum_{0\le\delta_1^\prime\le\delta_2^\prime} 
\omega(t(1,p^{\delta_1^\prime},p^{\delta_2^\prime})) p^{2\delta_1^\prime+\delta_2^\prime} (x_0{X})^{\delta_2^\prime}\\
&=[(1-x_0{X})(1-x_0x_1x_2x_3{X})]^{-1} \sum_{0\le\delta_1^\prime\le\delta_2^\prime} \omega(t(1,p^{\delta_1^\prime},p^{\delta_2^\prime})) 
p^{2\delta_1^\prime+\delta_2^\prime} (x_0{X})^{\delta_2^\prime}\\
&=[(1-x_0{X})(1-x_0x_1{X})(1-x_0x_2{X})(1-x_0x_3{X})(1-x_0x_1x_2{X})\\
&\ \ \ \ \cross (1-x_0x_1x_3{X})(1-x_0x_2x_3{X})(1-x_0x_1x_2x_3{X})]^{-1} P_3({X})
\end{split}
\end{equation}
Here $P_3({X})$ denotes a polynomial of degree 6 as stated in the Theorem 6 (page 451) of \cite{An70}.  
This rational polynomial presentation is proved in \cite{An69} for Hecke series and $\zeta$-functions
of the groups $GL_n$ and $SP_n$ over local fields.  Further theory and applications were developed 
for genus 2 in the work \cite{An74}.

It follows that
\begin{align}
P_3({X})&=
 \sum_{0\le\delta_1^\prime\le\delta_2^\prime} 
\omega(t(1,p^{\delta_1^\prime},p^{\delta_2^\prime})) p^{2\delta_1^\prime+\delta_2^\prime} (x_0{X})^{\delta_2^\prime}\times \\
&\nonumber \times [(1-x_0x_1{X})(1-x_0x_2{X})(1-x_0x_3{X})
\\ & \nonumber 
(1-x_0x_1x_2{X})(1-x_0x_1x_3{X})(1-x_0x_2x_3{X})],
\end{align}
and Theorem \ref{ThP} will follow from the explicit computation of  the coefficients 
$$\omega(t(1,p^{\delta_1^\prime},p^{\delta_2^\prime}))$$
given in the next section.

\section{
Images of the Hecke operators under the spherical map}
The formula for $P_3$ is obtained using the following computation for the images
$\omega ({t}(1,p^{\lambda_2},p^{\lambda_3}))$ of the Hecke operators under the spherical 
map for the group $\Lambda=GL_3(\Z)$.  
Note that the notation $\Omega$ used in the article \cite{An70}
corresponds to  $\omega$ in our formulas by the substitution of 
$x_1$ by $x_1/p$, $x_2$ by $x_2/p$ and $x_3$ by $x_3/p$.  
We used formula (1.7) for $\Omega$ at 
page 432 of \cite{An70} and adopted it for $\omega$.  
In the notations of that article we have that
$W$ is the group $S_n = S_3$, the set 
$$\Sigma = \{(1,-1,0),(1,0,-1),(0,1,-1)\}, q=p,\mbox{ and } 
\lambda = (0,\delta_1^{'},\delta_2^{'}).
$$  
The expression for the polynomial $c(x)$ from \cite{An70}
takes the following form
$$c(x_1,x_2,x_3)=\frac{(x_2-x_1/p)(x_3-x_1/p)(x_3-x_2/p)}{(x_2-x_1)(x_3-x_1)(x_3-x_2)}\,.$$
More precisely, let us define
$$
c_{\lambda, \mu}(x_1,x_2,x_3)=x_2^{\lambda}x_3^{\mu}\frac{(x_2-x_1/p)(x_3-x_1/p)(x_3-x_2/p)}{(x_2-x_1)(x_3-x_1)(x_3-x_2)}.
$$
Then $\ds
\omega ({t}(1,p^{\lambda},p^{\mu}))= K_{\lambda, \mu}\sum_{\sigma\in S_3}c_{\lambda, \mu}
(x_{\sigma (1)}, x_{\sigma (2)} ,x_{\sigma (3)})$, where
$$  K_{\lambda, \mu}=\frac 1{p^{2\lambda+\mu}}\times
\begin{cases} \ds \frac {p^3}{(p+1)(p^2+p+1)},& \mbox{ if }\lambda=\mu=0 \mbox{ (in this case } \omega ({t}(1,1,1))=1), \\
\ds \frac{p}{p+1}, & \mbox{ if }\lambda=0, \mu>0, 
 \mbox{ or }\lambda= \mu>0, \\
1,& \mbox{ otherwise }.
\end{cases}
$$

In the case of $n=3$ we have 28 different possibilities for $\omega(t(1,p^{\lambda_2},p^{\lambda_3}))$,
which we only need in order to compute the polynomial $P_3({X})$ of degree 6:

\begin{itemize}
\item[1)]
 $\omega(t(1,1,1))
 = 1$,
\item[2)]
 $\omega(t(1,1,p))
 = {\displaystyle \frac{\mathit{sym}_{1,0,0}}{p}}$,
\item[3)]
 $\omega(t(1,p,p))
 = {\displaystyle \frac{\mathit{sym}_{1,1,0}}{p^3}}$,
\item[4)]
 $\omega(t(1,1,p^2))
 = {\displaystyle \frac{(p-1)\,\mathit{sym}_{1,1,0}}{p^3}}
 + {\displaystyle \frac{\mathit{sym}_{2,0,0}}{p^2}}$,
\item[5)]
 $\omega(t(1,p,p^2))
 = {\displaystyle \frac{(2p^2-p-1)\,\mathit{sym}_{1,1,1}}{p^6}}
 + {\displaystyle \frac{\mathit{sym}_{2,1,0}}{p^4}}$,
\item[6)]
 $\omega(t(1,p^2,p^2))
 = {\displaystyle \frac{(p-1)\,{\mathit{sym}_{2,1,1}}}{p^7}}
 + {\displaystyle \frac{\mathit{sym}_{2,2,0}}{p^6}}$,
\item[7)]
 $\omega ({t}(1,1,p^{3}))
 = {\displaystyle \frac {(p^{2}  - 2\,p + 1)\,{\mathit{sym}_{1,1,1}}}{p^{5}}}
 + {\displaystyle \frac {(p^{2}-p)\,{\mathit{sym}_{2,1,0}}}{p^{5}}}
 + {\displaystyle \frac {{\mathit{sym}_{3,0,0}}}{p^{3}}}$,
\item[8)]
 $\omega ({t}(1,p,p^{3}))
 = {\displaystyle \frac {(2\,p-2)\,{\mathit{sym}_{2,1,1}}}{p^{6}}}
 + {\displaystyle \frac {(p - 1)\,{\mathit{sym}_{2,2,0}}}{p^{6}}}
 + {\displaystyle \frac {{\mathit{sym}_{3,1,0}}}{p^{5}}}$,
\item[9)]
 $\omega ({t}(1,p^{2},p^{3}))
 = {\displaystyle \frac {(2\,p-2)\,{\mathit{sym}_{2,2,1}}}{p^{8}}}
 + {\displaystyle \frac {(p - 1)\,{\mathit{sym}_{3,1,1}}}{p^{8}}} 
 + {\displaystyle \frac {{\mathit{sym}_{3,2,0}}}{p^{7}}}$, 
\item[10)]
 $\omega ({t}(1,p^{3},p^{3}))
 = {\displaystyle \frac {(p^{2} -2\,p + 1)\,{\mathit{sym}_{2,2,2}}}{p^{11}}} 
 + {\displaystyle \frac {(p^{2}-p)\,{\mathit{sym}_{3,2,1}}}{p^{11}}} 
 + {\displaystyle \frac {{\mathit{sym}_{3,3,0}}}{p^{9}}}$,
\item[11)]
 $\omega ({t}(1,1,p^{4}))
 = {\displaystyle \frac {(p^2 - 2\,p + 1)\,{\mathit{sym}_{2,1,1}}}{p^{6}}} 
 + {\displaystyle \frac {(p^{2}-p)\,{\mathit{sym}_{2,2,0}}}{p^{6}}} \\
 + {\displaystyle \frac {(p^{2}-p)\,{\mathit{sym}_{3,1,0}}}{p^{6}}} 
 + {\displaystyle \frac {{\mathit{sym}_{4,0,0}}}{p^{4}}}$,
\item[12)]
 $\omega ({t}(1,p,p^{4}))
\\
 = {\displaystyle \frac {(2\,p^{2} - 3\,p+1)\,{\mathit{sym}_{2,2,1}}}{p^{8}}} 
 + {\displaystyle \frac {(2\,p^{2} - 2\,p)\,{\mathit{sym}_{3,1,1}}}{p^{8}}} \\
 + {\displaystyle  \frac {(p^{2}-p)\,{\mathit{sym}_{3,2,0}}}{p^{8}}} 
 + {\displaystyle  \frac {{\mathit{sym}_{4, \,1, \,0}}}{p^{6}}}$,
\item[13)]
 $\omega ({t}(1,p^{2},p^{4}))
 = {\displaystyle \frac {( - 4\,p^{2} + 3\,p^{3} + 2\,p-1)\,{\mathit{sym}_{2,2,2}}}{p^{11}}} 
 + {\displaystyle \frac {(2\,p^3 - 3\,p^{2} + p)\,{\mathit{sym}_{3,2,1}}}{p^{11}}} \\	
 + {\displaystyle \frac {(p^{3} - p^{2})\,{\mathit{sym}_{3,3,0}}}{p^{11}}} 
 + {\displaystyle \frac {(p^{3} - p^{2})\,{\mathit{sym}_{4,1,1}}}{p^{11}}} 
 + {\displaystyle \frac {{\mathit{sym}_{4,2,0}}}{p^{8}}}$,
\item[14)]
 $\omega ({t}(1,p^{3},p^{4}))
 = {\displaystyle \frac {(2\,p^{2} - 3\,p+ 1)\,{\mathit{sym}_{3,2,2}}}{p^{12}}} 
\\ + {\displaystyle \frac {(2\,p^{2} - 2\,p)\,{\mathit{sym}_{3,3,1}}}{p^{12}}} 
 + {\displaystyle \frac {(p^{2}-p)\,{\mathit{sym}_{4,2,1}}}{p^{12}}} 
 + {\displaystyle \frac {{\mathit{sym}_{4,3,0}}}{p^{10}}}$,
\item[15)]
 $\omega ({t}(1,p^{4},p^{4}))
 = {\displaystyle \frac {(p^{2} -2\,p + 1)\,{\mathit{sym}_{3,3,2}}}{p^{14}}} 
 + {\displaystyle \frac {(p^{2}-p)\,{\mathit{sym}_{4,2,2}}}{p^{14}}} \\
 + {\displaystyle \frac {(p^{2}-p)\,{\mathit{sym}_{4,3,1}}}{p^{14}}} 
 + {\displaystyle \frac {{\mathit{sym}_{4,4,0}}}{p^{12}}}$,
\item[16)]
 $\omega (t(1,1,p^{5}))\\
 = {\displaystyle \frac {(p^{2} -2\,p+ 1)\,{\mathit{sym}_{2,2,1}}}{p^{7}}} 
 + {\displaystyle \frac {(p^{2} -2\,p+ 1)\,{\mathit{sym}_{3,1,1}}}{p^{7}}} 
 + {\displaystyle \frac {(p^{2}-p)\,{\mathit{sym}_{3,2,0}}}{p^{7}}} \\
 + {\displaystyle \frac {(p^{2}-p)\,{\mathit{sym}_{4,1,0}}}{p^{7}}} 
 + {\displaystyle \frac {{\mathit{sym}_{5,0,0}}}{p^{5}}}$,
\item[17)]
 $\omega ({t}(1,p,p^{5}))\\
 = {\displaystyle \frac {(2\,p^2 - 4\,p + 2)\,{\mathit{sym}_{2,2,2}}}{p^{9}}} 
 + {\displaystyle \frac {(2\,p^{2} - 3\,p+1)\,{\mathit{sym}_{3,2,1}}}{p^{9}}}  
 + {\displaystyle \frac {(p^{2}-p)\,{\mathit{sym}_{3,3,0}}}{p^{9}}} \\ 
 + {\displaystyle \frac {(2\,p^{2} - 2\,p)\,{\mathit{sym}_{4,1,1}}}{p^{9}}} 
 + {\displaystyle \frac {(p^{2}-p)\,{\mathit{sym}_{4,2,0}}}{p^{9}}} 
 + {\displaystyle \frac {{\mathit{sym}_{5,1,0}}}{p^{7}}}$,
\item[18)]
 $\omega ({t}(1,p^{2},p^{5}))\\
 = {\displaystyle \frac {(3\,p^3 - 5\,p^{2} + 3\,p - 1)\,{\mathit{sym}_{3,2,2}}}{p^{12}}} 
 + {\displaystyle \frac {(2\,p^3 - 4\,p^{2} + 2\,p)\,{\mathit{sym}_{3,3,1}}}{p^{12}}}  \\
 + {\displaystyle \frac {(2\,p^3 - 3\,p^{2} + p)\,{\mathit{sym}_{4,2,1}}}{p^{12}}} 
 + {\displaystyle \frac {(p^{3} - p^{2})\,{\mathit{sym}_{4,3,0}}}{p^{12}}} 
\\ + {\displaystyle \frac {(p^{3} - p^{2})\,{\mathit{sym}_{5,1,1}}}{p^{12}}}  
 + {\displaystyle \frac {{\mathit{sym}_{5,2,0}}}{p^{9}}}$,
\item[19)]
 $\omega ({t}(1,p^{3},p^{5}))\\
 = {\displaystyle \frac {(3\,p^3 - 5\,p^{2} + 3\,p - 1)\,{\mathit{sym}_{3,3,2}}}{p^{14}}} 
 + {\displaystyle \frac {(2\,p^3 - 4\,p^{2} + 2\,p)\,{\mathit{sym}_{4,2,2}}}{p^{14}}} \\  
 + {\displaystyle \frac {(2\,p^3 - 3\,p^{2} + p)\,{\mathit{sym}_{4,3,1}}}{p^{14}}} 
 + {\displaystyle \frac {(p^{3} - p^{2})\,{\mathit{sym}_{4,4,0}}}{p^{14}}} 
\\ + {\displaystyle \frac {(p^{3} - p^{2})\,{\mathit{sym}_{5,2,1}}}{p^{14}}} 
 + {\displaystyle \frac {{\mathit{sym}_{5,3,0}}}{p^{11}}}$,
\item[20)]
 $\omega ({t}(1,p^{4},p^{5}))
 = {\displaystyle \frac {(2\,p^2 - 4\,p + 2)\,{\mathit{sym}_{3,3,3}}}{p^{15}}} 
 + {\displaystyle \frac {(2\,p^2 - 3\,p + 1)\,{\mathit{sym}_{4,3,2}}}{p^{15}}} \\  
 + {\displaystyle \frac {(2\,p^2 - 2\,p    )\,{\mathit{sym}_{4,4,1}}}{p^{15}}} 
 + {\displaystyle \frac {(   p^2 -    p    )\,{\mathit{sym}_{5,2,2}}}{p^{15}}} 
 + {\displaystyle \frac {(   p^2 -    p    )\,{\mathit{sym}_{5,3,1}}}{p^{15}}} 
 + {\displaystyle \frac {                     {\mathit{sym}_{5,4,0}}}{p^{13}}}$,
\item[21)]
 $\omega ({t}(1,p^{5},p^{5}))
\\
 = {\displaystyle \frac {(p^2 - 2\,p + 1)\,{\mathit{sym}_{4,3,3}}}{p^{17}}} 
 + {\displaystyle \frac {(p^2 - 2\,p + 1)\,{\mathit{sym}_{4,4,2}}}{p^{17}}} 
 + {\displaystyle \frac {(p^2 -    p    )\,{\mathit{sym}_{5,3,2}}}{p^{17}}} \\ 
 + {\displaystyle \frac {(p^2 -    p    )\,{\mathit{sym}_{5,4,1}}}{p^{17}}} 
 + {\displaystyle \frac {                  {\mathit{sym}_{5,5,0}}}{p^{15}}}$,
\item[22)]
 $\omega ({t}(1,1,p^{6}))\\
 = {\displaystyle \frac {(p^2-2\,p+1)\,{\mathit{sym}_{2,2,2}}}{p^{8}}} 
 + {\displaystyle \frac {(p^2-2\,p+1)\,{\mathit{sym}_{3,2,1}}}{p^{8}}}  
 + {\displaystyle \frac {(p^2-   p  )\,{\mathit{sym}_{3,3,0}}}{p^{8}}} \\ 
 + {\displaystyle \frac {(p^2-2\,p+1)\,{\mathit{sym}_{4,1,1}}}{p^{8}}} 
 + {\displaystyle \frac {(p^2-   p  )\,{\mathit{sym}_{4,2,0}}}{p^{8}}}  
 + {\displaystyle \frac {(p^2-   p  )\,{\mathit{sym}_{5,1,0}}}{p^{8}}} 
 + {\displaystyle \frac {              {\mathit{sym}_{6,0,0}}}{p^{6}}}$,
\item[23)]
 $\omega ({t}(1,p,p^{6}))
 = {\displaystyle \frac {(2\,p^2 - 4\,p + 2)\,{\mathit{sym}_{3,2,2}}}{p^{10}}} 
 + {\displaystyle \frac {(2\,p^2 - 3\,p + 1)\,{\mathit{sym}_{3,3,1}}}{p^{10}}} \\ 
 + {\displaystyle \frac {(2\,p^2 - 3\,p + 1)\,{\mathit{sym}_{4,2,1}}}{p^{10}}} 
 + {\displaystyle \frac {(   p^2 -    p    )\,{\mathit{sym}_{4,3,0}}}{p^{10}}} 
 + {\displaystyle \frac {(2\,p^2 - 2\,p    )\,{\mathit{sym}_{5,1,1}}}{p^{10}}} \\
 + {\displaystyle \frac {(   p^2 -    p    )\,{\mathit{sym}_{5,2,0}}}{p^{10}}} 
 + {\displaystyle \frac {                     {\mathit{sym}_{6,1,0}}}{p^{8}}}$,
\item[24)]
 $\omega ({t}(1,p^{2},p^{6}))\\
 = {\displaystyle \frac {(3\,p^3 - 6\,p^2 + 4\,p - 1)\,{\mathit{sym}_{3,3,2}}}{p^{13}}}  
 + {\displaystyle \frac {(3\,p^3 - 5\,p^2 + 3\,p - 1)\,{\mathit{sym}_{4,2,2}}}{p^{13}}} \\
 + {\displaystyle \frac {(2\,p^3 - 4\,p^2 + 2\,p    )\,{\mathit{sym}_{4,3,1}}}{p^{13}}}  
 + {\displaystyle \frac {(   p^3 -    p^2           )\,{\mathit{sym}_{4,4,0}}}{p^{13}}} 
 + {\displaystyle \frac {(2\,p^3 - 3\,p^2 +    p    )\,{\mathit{sym}_{5,2,1}}}{p^{13}}} \\
 + {\displaystyle \frac {(   p^3 -    p^2           )\,{\mathit{sym}_{5,3,0}}}{p^{13}}} 
 + {\displaystyle \frac {(   p^3 -    p^2           )\,{\mathit{sym}_{6,1,1}}}{p^{13}}} 
 + {\displaystyle \frac {                              {\mathit{sym}_{6,2,0}}}{p^{10}}}$,
\item[25)]
 $\omega ({t}(1,p^{3},p^{6}))
 = {\displaystyle \frac {(4\,p^3 - 7\,p^2 + 5\,p - 2)\,{\mathit{sym}_{3,3,3}}}{p^{15}}}  
\\ + {\displaystyle \frac {(3\,p^3 - 6\,p^2 + 4\,p - 1)\,{\mathit{sym}_{4,3,2}}}{p^{15}}} 
 + {\displaystyle \frac {(2\,p^3 - 4\,p^2 + 2\,p    )\,{\mathit{sym}_{4,4,1}}}{p^{15}}}  
\\ + {\displaystyle \frac {(2\,p^3 - 4\,p^2 + 2\,p    )\,{\mathit{sym}_{5,2,2}}}{p^{15}}} 
 + {\displaystyle \frac {(2\,p^3 - 3\,p^2 +    p    )\,{\mathit{sym}_{5,3,1}}}{p^{15}}} 
 + {\displaystyle \frac {(   p^3 -    p^2           )\,{\mathit{sym}_{5,4,0}}}{p^{15}}}  
 + {\displaystyle \frac {(   p^3 -    p^2           )\,{\mathit{sym}_{6,2,1}}}{p^{15}}} 
 + {\displaystyle \frac {                              {\mathit{sym}_{6,3,0}}}{p^{12}}}$,
\item[26)]
 $\omega ({t}(1,p^{4},p^{6}))\\
 = {\displaystyle \frac {(3\,p^3 - 6\,p^2 + 4\,p - 1)\,{\mathit{sym}_{4,3,3}}}{p^{17}}}  
 + {\displaystyle \frac {(3\,p^3 - 5\,p^2 + 3\,p - 1)\,{\mathit{sym}_{4,4,2}}}{p^{17}}} \\
 + {\displaystyle \frac {(2\,p^3 - 4\,p^2 + 2\,p    )\,{\mathit{sym}_{5,3,2}}}{p^{17}}}  
 + {\displaystyle \frac {(2\,p^3 - 3\,p^2 +    p    )\,{\mathit{sym}_{5,4,1}}}{p^{17}}} 
 + {\displaystyle \frac {(   p^3 -    p^2           )\,{\mathit{sym}_{5,5,0}}}{p^{17}}} \\
 + {\displaystyle \frac {(   p^3 -    p^2           )\,{\mathit{sym}_{6,2,2}}}{p^{17}}} 
 + {\displaystyle \frac {(   p^3 -    p^2           )\,{\mathit{sym}_{6,3,1}}}{p^{17}}} 
 + {\displaystyle \frac {                              {\mathit{sym}_{6,4,0}}}{p^{14}}}$,
\item[27)]
 $\omega ({t}(1,p^{5},p^{6}))
 = {\displaystyle \frac {(2\,p^2 - 4\,p + 2)\,{\mathit{sym}_{4,4,3}}}{p^{18}}} 
 + {\displaystyle \frac {(2\,p^2 - 3\,p + 1)\,{\mathit{sym}_{5,3,3}}}{p^{18}}} \\ 
 + {\displaystyle \frac {(2\,p^2 - 3\,p + 1)\,{\mathit{sym}_{5,4,2}}}{p^{18}}} 
 + {\displaystyle \frac {(2\,p^2 - 2\,p    )\,{\mathit{sym}_{5,5,1}}}{p^{18}}} 
 + {\displaystyle \frac {(   p^2 -    p    )\,{\mathit{sym}_{6,3,2}}}{p^{18}}} \\
 + {\displaystyle \frac {(   p^2 -    p    )\,{\mathit{sym}_{6,4,1}}}{p^{18}}} 
 + {\displaystyle \frac {                     {\mathit{sym}_{6,5,0}}}{p^{16}}}$,
\item[28)]
 $\omega ({t}(1,p^{6},p^{6}))
 = {\displaystyle \frac {(p^2 - 2\,p + 1)\,{\mathit{sym}_{4,4,4}}}{p^{20}}} 
 + {\displaystyle \frac {(p^2 - 2\,p + 1)\,{\mathit{sym}_{5,4,3}}}{p^{20}}} \\
 + {\displaystyle \frac {(p^2 - 2\,p + 1)\,{\mathit{sym}_{5,5,2}}}{p^{20}}} 
 + {\displaystyle \frac {(p^2 -    p    )\,{\mathit{sym}_{6,3,3}}}{p^{20}}} 
 + {\displaystyle \frac {(p^2 -    p    )\,{\mathit{sym}_{6,4,2}}}{p^{20}}} \\
 + {\displaystyle \frac {(p^2 -    p    )\,{\mathit{sym}_{6,5,1}}}{p^{20}}} 
 + {\displaystyle \frac {                  {\mathit{sym}_{6,6,0}}}{p^{18}}}$ \, . \qed
\end{itemize}

Note that Rhodes, J. A. and Shemanske, T. R. developed an alternative method of computing
$\omega(t(p^{\delta_1},\dots ,p^{\delta_n}))$ in \cite{RhSh}, based on counting of certain left
cosets in a given double coset (Theorem 4.3 of \cite{RhSh}). 

\end{proof}

\section{A special case}
For some particular values of the Satake parameters $x_0$, $x_1$, $x_2$, $x_3$, the polynomial $P_3$
can be considerably simplified.  For example,
let us substitute $x_0=1$, $x_1=p$, $x_2=p^2$ and $x_3=p^3$ as in Exercise 3.3.40, p.168 of \cite{An87}:
$\mathit{P}_\nu ({X}) :=\mathit{P}(1, p, p^2, p^3, {X})$, where $\nu$ denotes the degree homomorphism
$\nu(x_0)=1$, $\nu(x_1)=p$, $\nu(x_2)=p^2$, $\nu(x_3)=p^3$. Then the polynomial $P$ takes the form

\begin{align*}&
\mathit{P}_\nu ({X}) = 1 - \left({\displaystyle \frac {p^{7} + p^{8} + p^{9
}}{p}} + (p^{2} + p + 1)\,p^{4} + {\displaystyle \frac {p^{3} + 
p^{4} + p^{5}}{p}} \right)\,{X}^{2} 
\\ &
\mbox{} + \left((p + 1)\,p^{10} + {\displaystyle \frac {(p + 1)\,(p^{9
} + p^{10} + p^{11})}{p^{2}}}  + {\displaystyle \frac {(p + 1)\,(
p^{7} + p^{8} + p^{9})}{p^{2}}}  + (p + 1)\,p^{4}\right)\,{X}^{3} \\ &
\mbox{} - \left({\displaystyle \frac {p^{13} + p^{14} + p^{15}}{p^{
2}}} + (p^{2} + p + 1)\,p^{9}+ {\displaystyle \frac {p^{9} + p
^{10} + p^{11}}{p^{2}}} \right)\,{X}^{4} + p^{15}\,{X}^{6} \\ &
=1 - (p^{8}\,+ p^{7}\,+2\,p^{6}\,+ p^{5}\,+2\,p^{4}\,+ p^{3}+p^{2}){X}^{2}  \\ &
+(p^{11}+ 2\,p^{10} + 2\,p^{9}+ 3\,p^{8}   + 3\,p^{7} + 2\,p^{6} + 
2\,p^{5} + p^{4})\,{X}^{3} \\ &
-(p^{13}\,+ p^{12}\,+2\,p^{11}\,+ p^{10}\,+2\,p^{9}\,+ p^{8}+p^{7}){X}^{4}
+ p^{15}\,{X}
^{6}. 
\end{align*}
This gives the following factorization:
\begin{align*}&
\mathit{P}_\nu ({X})
=(1-p^{}\,{X})\,(1-p^{2}\,{X} )\,(1-p^{3}\,{X} )\,(1-p^{4}\,{X})\,
\\ & \hskip3cm
\cross(1+ p^{}\,{X} + p^{2}\,{X} + p^{3}\,{X} + p^{4}\,{X} +p^{5}\,{X}^{2} ).
\end{align*}

\subsection*{Acknowledgements}
We are very grateful to the Institute Fourier
 (UJF, Grenoble-1) 
for  the  permanent excellent working environment. 

It is a great pleasure for us to thank Siegfried Boecherer, Gilles Robert and Tanguy Rivoal 
 for valuable  discussions and observations. 

Our results were presented in  a 
 talk   \cite{PaGRFA} on June 22, 2006,
``Produits d'Euler attach\'es aux formes modulaires
de Siegel et la conjecture de Shimura explicite  pour $Sp_3$''
at  S\'eminaire 
Groupes R\'eductifs et Formes Automorphes
\`a l'Institut de
Math\'ematiques de Jussieu,
(Chevaleret) 
 by the first author, thankful to Michael Harris for the invitation.

Our special thanks go to Anatoli Nikolayevich Andrianov and Yuri Ivanovich Manin, 
 for providing us with  advice and encouragement.


\bibliographystyle{url}

\end{document}